\documentclass{amsart}[12pt]
\usepackage{amsmath, amsthm, amscd, amsfonts}

\begin{document}

\title[CONVEXITY OF \u{C}EBY\u{S}EV SETS]{CONVEXITY OF \u{C}EBY\u{S}EV SETS THROUGH DIFFERENTIABILITY OF DISTANCE FUNCTION }

\maketitle

\author \hspace{0.75cm} Amanollah Assadi, Hadi Haghshenas and
Hossein Hosseini Guive

\begin{center}
Department of Mathematical Analysis, Birjand University, Birjand, Iran\\
aman$_{-}$assadi@yahoo.com\\haghshenas60@gmail.com\\hossein.giv@gmail.com

\end{center}

\begin{abstract}
In this paper, we study a part of approximation theory that
presents the conditions under which a \u{C}eby\u{s}ev set in a
Banach space is convex. To do so, we use Gateaux differentiability
of the distance function.
\end{abstract}


\textbf{\textit{Key Words:}}\hspace{0.5cm} Distance function,
nearest point, \u{C}eby\u{s}ev set, strictly convex space, smooth
space, Gateaux differentiability, subdifferential.

\section{Introduction}
The approximation theory is one of the important branch of
functional analysis that \u{C}eby\u{s}ev originated it in
nineteenth century. But, convexity of \u{C}eby\u{s}ev sets is one
of the basic problems in this theory. In a finite dimensional
smooth normed space a \u{C}eby\u{s}ev set is convex and for
infinite dimensional, every weakly closed \u{C}eby\u{s}ev set in a
smooth and uniformly convex Banach space is convex. Every
boundedly compact \u{C}eby\u{s}ev set in a smooth Banach space is
convex and in a Banach space, which is uniformly smooth, each
approximately compact \u{C}eby\u{s}ev set is convex (The concept
of approximatively compact sets introduced by N. V. Efimov and S.
B. Stechkin), and that in a strongly smooth space or in a Banach
space $X$ with strictly convex dual $X^{*}$, every \u{C}eby\u{s}ev
set with continuous metric projection is convex, (\cite{2}). There
are still several open problems concerning convexity of
\u{C}eby\u{s}ev sets. Can we prove that in some Banach spaces, a
nonempty subset is a \u{C}eby\u{s}ev set if and only if it is
closed and convex? This is unsolved, even in the special case of
infinite-dimensional Hilbert space. In 1934, L. N. H. Bunt proved
that each \u{C}eby\u{s}ev set in a finite-dimensional Hilbert
space must be convex. From this result, we see that in a
finite-dimensional Hilbert space, a nonempty subset is a
\u{C}eby\u{s}ev set if and only if it is closed and convex. In
\cite{5}, G. G. Johnson gave an example: there exists an
incomplete inner product space which possesses a non-convex
\u{C}eby\u{s}ev set (M. Jiang completed the proof in 1993). Is
there an infinite-dimensional Hilbert space possessing a
non-convex \u{C}eby\u{s}ev set? As addressed above, it is unknown.
In the last part of the paper, we present some
conditions under which a \u{C}eby\u{s}ev subset is convex.\\

\section{Basic definitions and Preliminaries}
In this section, we collect some elementary facts which will help
us to establish our main results. For details the reader is
referred to \cite{4}. As the first step, let us fix our notation.
Through this paper, $(X,\|.\|)$ denotes a real Banach space and
$S(X)=\{x\in X;\|x\|=1\}$.\\For an element $x \in X$ and a
nonempty subset $K$ in X, we define the distance function $d_{K}:
X \rightarrow \Bbb{R}$ by $d_{K}(x)= inf \{\|y-x\| ; y\in K \}$.
It is easy to see that the value of $d_{K}(x)$ is zero if and only
if $x$ belongs to $\overline{K}$, the closure of $K$. The subset
$K$ is called proximinal (resp. \u{C}eby\u{s}ev), if for each $x
\in X \setminus K$, the set of best approximations to $x$ from $K$
$$P_{K}(x)= \{y \in K; \|y-x\|=d_{K}(x)\},$$is nonempty (resp. a singleton). This concept was introduced by S. B. Stechkin
and named after the founder of best approximation theory,
\u{C}eby\u{s}ev.\\One interesting and fruitful line of research,
dating from the early days of Banach space theory, has been to
relate analytic properties of a Banach space to various
geometrical conditions on the Banach space. The simplest example
of such a condition is that of strict convexity. It is often
convenient to know whether the triangle inequality is strict for
non collinear points in a given Banach space. We say that the norm
$\|.\|$ of $X$ is strictly convex (rotund) if,
$$\|x+y\|<\|x\|+\|y\|$$whenever $x$ and $y$ are not parallel. That
is, when they are not multiples of one another. \\Related to the
notion of strict convexity, is the notion of smoothness.\\ We say
that, the norm $\|.\|$ of $X$ is smooth at $x \in X \setminus
\{0\}$ if, there is a unique $f \in X^{*}$ such that $\|f\|=1$ and
$f(x)=\|x\|$. Of course, the Hahn-Banach theorem ensures the
existence of at least one such functional $f$.\\The spaces
$L^{p}(\mu), 1<p<\infty$, are strictly convex and smooth, while
the spaces $L^{1}(\mu)$ and $C(K)$ are neither strictly convex nor
smooth except in the trivial case when they are one
dimensional.\\If the dual norm of $X^{*}$ is smooth, then the norm
of $X$ is strictly convex and if the dual norm of $X^{*}$ is
strictly convex, then the norm of $X$ is smooth. Note that, The
converse is true only for reflexive spaces. There are examples of
strictly convex spaces whose duals fail to be smooth.\\Let
$f:X\rightarrow \Bbb{R}$ be a function and $x,y\in X$. Then $f$ is
said to be Gateaux differentiable at $x$ if, there exists a
functional $A\in X^{*}$ such that
$A(y)=\displaystyle{\lim_{t\rightarrow 0}\frac{f(x+ty)-f(x)}{t}}$.
In this case $f$ is called Gateaux differentiable at $x$ with the
Gateaux derivative $A$ and $A$ is denoted by $f'(x)$. In this
case, the $A(y)$ is denoted usually by $<f'(x),y>$. If the limit
above exists uniformly for each $y\in S(X)$, then $f$ is
Fr\'{e}chet differentiable at $x$ with Fr\'{e}chet derivative $A$.
Similarly, the norm function $\|.\|$ is Gateaux (Fr\'{e}chet)
differentiable at non-zero $x$ if the function $f(x)=\|x\|$ is
Gateaux differentiable.\\In the general, Gateaux differentiability
not imply Fr\'{e}chet differentiability. For example the canonical
norm of $l^{1}$ is nowhere Fr\'{e}chet differentiable and it is
Gateaux differentiable at $x=(x_{i})_{i \in \Bbb{N}}$ if and only
if $x_{i}\neq 0 $ for every $i \in \Bbb{N}$.\\The norm of any
Hilbert space, is Fr\'{e}chet differentiable at nonzero points.\\
Suppose $f:X\rightarrow \Bbb{R}$ is a function and $x\in X$. The
functional $x^{*}\in X^{*}$ is called a subdifferential of $f$ at
$x$ if $\langle x^{*},y-x \rangle\leq f(y)-f(x)$, for all $y\in
X$. The set of all subdifferentials of $f$ at $x$ is denoted by
$\partial f(x)$ and we say that $f$ is
subdifferentiable at $x$ if $\partial f(x)\neq \emptyset$.\\
The following theorems presents relationship between various
notions of differentiability for norm and
the properties of the related space. \\
\textbf{Theorem 1. }\cite{4} The norm $\|.\|$ is Gateaux
differentiable at $x\in X \setminus \{0\}$ if and only if $X$ is
smooth in $x$.\\\textbf{Theorem 2. }\cite{4} If the dual norm of
$X^{*}$ is Fr\'{e}chet differentiable, then $X$ is reflexive.\\
\textbf{Theorem 3. }\cite{4} Let $f:X\rightarrow \Bbb{R}$ be a
convex function continuous at $x\in X$ and $\partial f(x)$ is a
singleton. Then $f$ is Gateaux differentiable at $x$.\\In the last
theorem, notice that continuity of $f$ in $x$ is an essential
condition. For example, if $f(x)=1+sin\frac{1}{x}$ for all $ x
\neq 0$ and $f(0)=0$, then $f$ is not continuous at $x=0$. Also
$\partial
f(0)=\{0\}$, while $f$ is not Gateaux differentiable at $x=0$.\\
For a real-valued function $\phi$ on $X$ and $x\in X$, set
$$F_{\phi}(x)=\sup_{\|y\|=1}\sup_{z\in X}\limsup_{t\rightarrow0^{+}}
\frac{\phi(x+tz+ty)-\phi(x+tz)}{t}.$$\\ \textbf{Lemma 1. }\cite{3}
Let $\phi$ is a real-valued function on $X$, $x\in X$ and
$y_{0}\in S(X)$ such that the Gateaux derivative of $\phi$ in $x$
exists and $\langle \phi'(x),y_{0} \rangle =F_{\phi}(x)$. If the
norm of $X$ is Gateaux differentiable at $y_{0}$ with Gateaux
derivative $f_{y_{0}}$, then $\phi$ is Gateaux differentiable at
$x$ and for each $y\in X$ we have $\langle \phi'(x),y \rangle
=F_{\phi}(x)f_{y_{0}}(y)$.\\\\Now the Lemma 1, give us the
following corollary, since distance functions are Lipschitz.:\\
\textbf{Corollary 1. }Let $K\subseteq X$ is closed and non-empty,
$x\in X\backslash K$ and $\overline{x}$ is a nearest point for $x$
in $K$, the distance function $d_{K}$ is Gateaux differentiable at
$x$ and in the direction of $(x-\overline{x})$ and the norm of $X$
is Gateaux differentiable at $(x-\overline{x})$ with Gateaux
derivative $f_{(x-\overline{x})}$. Then $d_{K}$ is Gateaux
differentiable at $x$ and $\langle d'_{K}(x),y \rangle =
f_{(x-\overline{x})}(y)$ for all $y\in
X$.\\
For nonempty closed subset $K$ of $X$ and $x,y\in X$, set
$$d_{K}^{-}(x;y)=\liminf_{t\rightarrow0^{+}}\frac{d_{K}(x+ty)-d_{K}(x)}{t}$$\\and\\
$$d_{K}^{+}(x;y)=\limsup_{t\rightarrow0^{+}}\frac{d_{K}(x+ty)-d_{K}(x)}{t}.$$
\\\textbf{Corollary 2.   }Suppose $K\subseteq X$ is closed and
nonempty, $x\in X\backslash K$, $\overline{x}$ is a nearest point
for $x$ in $K$. If the norm of $X$ is Gateaux differentiable at
$(x-\overline{x})$ and $d_{K}^{-}(x;x-\overline{x})=d_{K}(x)$,
then $d_{K}$ is Gateaux
differentiable at $x$.\\
\textit{Proof.   } Due to the norm of $X$ is Gateaux
differentiable at $x-\overline{x}$, it is sufficient to prove the
existence of the limit
$$\displaystyle{\lim_{t\rightarrow0}\frac{d_{K}(x+t(x-\overline{x}))-d_{K}(x)}{t}}.$$
Since $d'_{K}(x;\overline{x}-x)=-d_{K}(x)$ and
$\displaystyle{\lim_{t\rightarrow0^{-}}\frac{d_{K}(x+t(x-\overline{x}))-d_{K}(x)}{t}=d_{K}(x)},$
it is sufficient to prove that
$$\displaystyle{\lim_{t\rightarrow0^{+}}\frac{d_{K}(x+t(x-\overline{x}))
-d_{K}(x)}{t}=d_{K}(x)}.$$ For each $t>0$ we have
$d_{K}(x+t(x-\overline{x}))-d_{K}(x)\leq td_{K}(x).$ Hence,
$d_{K}^{+}(x;x-\overline{x})\leq d_{K}(x)$. If
$d_{K}^{-}(x;x-\overline{x})= d_{K}(x)$, then
$d_{K}(x)=d_{K}^{-}(x;x-\overline{x})\leq
d_{K}^{+}(x;x-\overline{x})\leq d_{K}(x).$ It follows
that $< d'_{K}(x),x-\overline{x}>$ exists and is equal to $d_{K}(x)$.\\
\\\textbf{Theorem 4.   }\cite{4} If the dual space of $X$ is strictly convex, then each closed nonempty subset $K$ in $X$ satisfying
$\displaystyle{\limsup_{\|y\| \rightarrow
0}\frac{d_{K}(x+y)-d_{K}(x)}{\|y\|}=1}$ for all $(x\in X\backslash
K)$ is convex.\section{Main Result}
 We start this section with our main result.\\
 \textbf{Theorem 5. }Suppose the dual space of $X$ is strictly convex, $K\subseteq X$ is
a \u{C}eby\u{s}ev set, $x\in X\backslash K$ and $\partial d_{K}(x)$ is singleton. The following are equivalent:\\
(\textbf{i}) $K$ is convex.\\ (\textbf{ii}) $d_{K}$ is convex .\\ (\textbf{iii}) $d_{K}$ is Gateaux differentiable at $x$.\\
(\textbf{iv}) There is $z\in S(X)$ such that $\displaystyle{\lim_{t\rightarrow 0^{+}}\frac{d_{K}(x+tz)-d_{K}(x)}{t}=1}$.\\
(\textbf{v}) $\displaystyle{\limsup_{\|y\| \rightarrow
0}\frac{d_{K}(x+y)-d_{K}(x)}{\|y\|}=1}$.\\\\\textit{Proof.   }
$(\textbf{i} \Rightarrow\textbf{ii}) $ Since $K$ is closed convex
set, $d_{K}$ is convex\cite{4}. \\ $(\textbf{ii} \Rightarrow
\textbf{iii})$ Since $d_{K}$ is convex and continuous at $x$ and
$\partial d_{K}(x)$ is singleton, $d_{K}$ is Gateaux
differentiable at $x$ and $\{d'_{K}(x)\}=\partial d_{K}(x)$.\\
$(\textbf{iii} \Rightarrow \textbf{iv})$ First note that by the
definition of \u{C}eby\u{s}ev sets there is a unique element
$\overline{x}\in K$ such that $\|x-\overline{x}\|=d_{K}(x)$. It
follows from Gateaux differentiability of $d_{K}$ that,
$\displaystyle{\liminf_{t\rightarrow
0^{+}}\frac{d_{K}(x+ty)-d_{K}(x)}{t}}$ exists for every $y \in X$.
For each $t>0$ we have, $d_{K}(x+t(x-\overline{x}))-d_{K}(x)\leq
td_{K}(x).$ Hence for $y=x-\overline{x}$, we
set:$$\displaystyle{\liminf_{t\rightarrow
0^{+}}\frac{d_{K}(x+t(x-\overline{x}))-d_{K}(x)}{t}=d_{K}(x)}.$$
Since $x\in X\backslash K$, $d_{K}(x)>0$ and if
$\displaystyle{t'=\frac{t}{d_{K}(x)}}$ as $t\rightarrow 0^{+}$,
Then by the above: $$\displaystyle{\liminf_{t'\rightarrow
0^{+}}\frac{d_{K}(x+t'(x-\overline{x}))-d_{K}(x)}{t'}=d_{K}(x)},$$
If now
$\displaystyle{z=\frac{x-\overline{x}}{\|x-\overline{x}\|}}$, then
$\|z\|=1$ and we have $\displaystyle{\liminf_{t\rightarrow
0^{+}}\frac{d_{K}(x+tz)-d_{K}(x)}{t}=1}$. On the other hand,
$d_{K}$ is a Lipschitz function and so
$\displaystyle{\limsup_{t\rightarrow
0^{+}}\frac{d_{K}(x+tz)-d_{K}(x)}{t}\leq 1}$.\\ $(\textbf{iv}
\Rightarrow \textbf{v})$ Since $d_{K}$ is a Lipschitz function
$\displaystyle{\limsup_{\|y\|\rightarrow
0}\frac{d_{K}(x+y)-d_{K}(x)}{\|y\|}\leq 1}.$ On the other hand for
each $v\in S(X)$,
$$\lim_{t \rightarrow 0^{+}}\frac{d_{K}(x+tv)-d_{K}(x)}{t}\leq\limsup_{\|y\|\rightarrow 0}\frac{d_{K}(x+y)-d_{K}(x)}{\|y\|} ,$$ \\
in particular for $v=z$ in \textbf{iv}, we have $1\leq
\displaystyle{\limsup_{\|y\|\rightarrow
0}\frac{d_{K}(x+y)-d_{K}(x)}{\|y\|}} .$\\
$(\textbf{v} \Rightarrow \textbf{i})$ This follows from
\textbf{theorem 4}.\\
\textbf{Remark 1.   }Suppose that the norm of $X$ and the dual
norm of $X^{*}$ are Fr\'{e}chet differentiable, $K\subseteq X$ is
\u{C}eby\u{s}ev and $x\in X\backslash K$. Then $X$ is reflexive,
since the dual norm of $X^{*}$ is Fr\'{e}chet differentiable.
Moreover $X$ is smooth, since the norm of $X$ is Fr\'{e}chet
differentiable. Thus $X^{*}$ is strictly convex. If now $d_{K}$ is
Gateaux differentiable at $x$, then $K$ is
convex.\\\\
\textbf{Remark 2.   } Suppose that $K\subseteq X$ is
\u{C}eby\u{s}ev, $x\in X\backslash K$ and $X^{*}$ is strictly
convex. By the definition of \u{C}eby\u{s}ev sets, there is unique
$\overline{x}\in K$ such that $\|x-\overline{x}\|=d_{K}(x)$. If
now $d_{K}^{-}(x;x-\overline{x})=d_{K}(x)$, then by
\textbf{corollary 2 } and \textbf{Remark 1}, $K$ is convex.\\

\end{document}